\documentclass[ijoc,nonblindrev]{informs3} 

\OneAndAHalfSpacedXII 

\usepackage{float}
\usepackage{geometry}
\usepackage[none]{hyphenat}
\usepackage{verbatim}
\usepackage[ruled,linesnumbered]{algorithm2e}
\usepackage{scrextend}
\usepackage{mathtools}
\usepackage[table,xcdraw]{xcolor}

\def \alga {$\mathcal{A}$}
\def \algb {$\mathcal{B}$}
\def \algc {$\mathcal{C}$}
\def \algd {$\mathcal{D}$}
\def \alge {$\mathcal{E}$}

\usepackage{natbib}
 \def\bibfont{\small}%
\TheoremsNumberedThrough     

\EquationsNumberedThrough    

\MANUSCRIPTNO{JOC-2021-02-OA-049} 

\begin{document}


\RUNAUTHOR{Baxter, Keskinocak, and Singh}

\RUNTITLE{Heterogeneous Multi-Resource Allocation with Subset Demand Requests}

\TITLE{Heterogeneous Multi-Resource Allocation with Subset Demand Requests}

\ARTICLEAUTHORS{%
\AUTHOR{Arden Baxter$^{1,2}$, Pinar Keskinocak$^{1,2}$, Mohit Singh$^{1}$} 
\AFF{$^1$H. Milton Stewart School of Industrial and Systems Engineering and $^2$ Center for Health and Humanitarian Systems, Georgia Institute of Technology, Atlanta, GA, 30332,
    \EMAIL{abaxter@gatech.edu}}
} 

\ABSTRACT{%
We consider the problem of allocating multiple heterogeneous resources geographically and over time to meet demands that require some subset of the available resource types simultaneously at a specified time, location, and duration. The objective is to maximize the total reward accrued from meeting (a subset of) demands. We model this problem as an integer program, show that it is NP-hard, and analyze the complexity of various special cases. We introduce approximation algorithms and an extension to our problem that considers travel costs. Finally, we test the performance of the integer programming model in an extensive computational study.  
}%



\maketitle

%


\section{Introduction}\label{intro}
In this paper, we introduce the heterogeneous multiple resource type allocation problem where multiple (e.g., a subset) resource types are requested by demands ($mRmD$) simultaneously at a specified time and location for a certain duration and the goal is to maximize total reward from meeting (as subset of) demands. $mRmD$ has vast applications in resource allocation and scheduling. For example, hospital operations require the coordinated scheduling of doctors, nurses, and operating rooms.
Similarly, in the home health care setting (i.e., providing supportive care in the home for illness, injury, or disability), patients may require visits from various members of the home health care team (e.g., home health aides, registered nurses, therapists, physicians, etc.). In some cases, multiple members may be needed simultaneously; for example, home health aides may need direct supervision by a registered nurse to perform any task for which they have not received satisfactory training (\cite{hha}).

Resource allocation and scheduling in networks have been studied across various applications in the operations research literature, including vehicle routing, machine scheduling, and robotics task allocation. However, most of the previous work considered either a single resource type (e.g., vehicle, ambulance, commodity, etc.) (\cite{deangelis,huang}) or independently scheduling/routing multiple resource types (\cite{viswanath}). There remains a significant gap in the literature on the efficient allocation of resources that may require some level of collaboration or coordination when the requirements of a certain demand cannot be met by a single resource type, i.e., resources of distinct types may be needed simultaneously or sequentially to meet a demand. While some studies (\cite{rauchecker,altay,Su,Lee}) explored the idea of collaboration between resources, the majority focused on tasks that can be done sequentially or independently by different resources. 


 We model $mRmD$ as an integer program, show that the problem and some special cases are NP-hard while other special cases are solvable in polynomial time, present approximation algorithms with provable bounds, and present results of a computational study which builds on the theoretical foundation. 

The remainder of this paper is structured as follows: Section \ref{lit_review} presents and discusses relevant literature. Sections \ref{prob_desc} and \ref{formulation_section} introduce our problem and give a description of the formulation. Theoretical results are presented in Sections \ref{complexity_results} and \ref{algorithms}. In Section \ref{prob_ext}, we introduce an extension to our problem that considers travel costs. Finally, a computational study is discussed in Section \ref{comp_study} before concluding in Section \ref{conclusion}.

\section{Literature Review}\label{lit_review}
 Problems  similar to $mRmD$ have been studied in emergency response and disaster management (\cite{baxter}), machine and project scheduling, robotic task allocation, and vehicle routing problems. Some of the theoretical notions of our work also align with interval graphs (\cite{gilmore,stephan,yannakakis,mertzios,carlisle}) where a demand can be represented as an interval with its start time and service duration. 

Multi-task scheduling aims to schedule jobs with different tasks among multiple machines to minimize the maximum completion time (makespan) (\cite{weizhen}) or maximize utility (\cite{fang}). These problems typically do not consider a spatial component (or sequence-dependent setup times) nor multiple simultaneous resource requirements. An exception is \cite{chen}, studying the one-job-on-multiple-machines model where several machines need to be assigned  simultaneously to process each job to minimize the completion time of all jobs.  

Within the emergency/humanitarian response management and vehicle routing streams, researchers addressed the problem of scheduling/routing single or multiple heterogeneous resources, with some coordination or collaboration in the latter case (\cite{rauchecker,altay,Su,Lee}). \cite{deangelis,viswanath,huang} studied vehicle routing and resource allocation decisions during humanitarian relief, with a single resource type and no dependencies among multiple resource types. \cite{rauchecker} defined collaboration as \textit{tight} (all resource units needed by a demand must arrive simultaneously) or \textit{loose} (resources may work independently to meet demand). While most of the previous work focused on loose collaboration, $mRmD$ calls for tight collaboration.  \cite{doulabi,bredstrom,mascolo} focused on the vehicle routing problem with synchronized visits (an example of tight collaboration where customers may require multiple vehicles simultaneously) and considered objectives including minimizing waiting time/delay in service, travel time, or costs.

The multi-robot task allocation problem is to allocate several homogeneous or heterogeneous robots (resources) to a number of tasks (demands) under system constraints to minimize the makespan (\cite{xu,gombolay,zheng,liu}), maximize the utility (\cite{nelke}), or minimize the total distance traveled (\cite{kartal}). \cite{xu,zheng,liu} characterize tasks by spatial constraints whereas \cite{nelke,gombolay,kartal} consider both spatial and temporal constraints.  \cite{xu,nelke,zheng} consider multi-robot tasks, which may require multiple homogeneous (\cite{nelke,xu,zheng}) or heterogeneous (\cite{nelke}) robots to be completed. In contrast to our study, in (\cite{nelke}) a robot can delay the start or interrupt the service of a task for a penalty, and all tasks are eventually completed, with the goal of maximizing the total utility.

\section{Problem Description}\label{prob_desc}
 A set of demands ($D$)  has to be processed by a set of resource types ($R$). Each resource type $r \in R$ has a set of starting locations, $S^r$, and $l_s^r$ resources at each starting location $s \in S^r$, i.e., there are  $\sum\limits_{s \in S^r} l_s^r$ resources of type $r \in R$. A demand $d \in D$ requires a single unit of resource types $M_d \subseteq R$ simultaneously at time $\tau_d$ for a duration of $\Delta_d$ (service time), resulting in a reward of $w_d$ only if met on time. Travel time between any two locations $x$ and $y$ is denoted by $f_{xy}$ (e.g., $f_{ij}$ denotes the travel time between the locations of demands  $i$ and $j$ and $f_{sd}$ denotes the travel time between a resource's starting location $s$ and the location of demand $d$). The goal is to assign resources to demands to maximize the total reward of meeting (a subset of) demands. The notation is summarized in Table~\ref{formulation}.
 
 \begin{definition}
 An instance $\mathcal{I}$ of $mRmD$ is defined as (i) sets $D,R,S^r$ for all $r \in R$, and $M_d$ for all $d \in D$, (ii) integer-valued parameters $l_s^r$ for all $s \in S^r, r \in R$ and $\tau_d,\Delta_d,w_d$ for all $d \in D$, and (iii) integer-valued travel times, $f_{xy}$, between any two locations $x$ and $y$.
 \end{definition}

\textbf{Pre-processing:} Note that an instance $\mathcal{I}$ of $mRmD$ may be described as a directed acyclic graph. Nodes represent the resource starting locations $S^r, r \in R$ and demands $d \in D$. Then, we may create arcs $(i,j)$ between two demand nodes $i$ and $j$ if and only if demand $j$ can be served immediately after demand $i$ by a resource, i.e., $\tau_i + \Delta_i + f_{ij} \leq \tau_j$. Further, we create arcs $(s,d)$ between a resource starting location node $s$ and a demand node $d$ if and only if a resource from that starting location can serve demand $d$, i.e.,  $f_{sd} \leq \tau_d$. We introduce pre-processing adjacency matrices $A$ and $B$ to determine the arcs $(i,j)$ and $(s,d)$ created in the directed graph, respectively. That is, matrix $A$ has a row and column for each demand and cell $A_{ij} = 1$ if and only if demand  $j$ can be served immediately after demand $i$ and the resource type requirements of $i$ and $j$ intersect. Matrix $B$ has a row for each resource's starting location,  a column for each demand, and cell $B_{sd} = 1$ if and only if a resource at starting at location $s$ can serve demand $d$ (on time, if $d$ was the first demand to be served by this resource). The pre-processing adjacency matrices are formally defined as follows: \begin{align*}
    A_{ij} &= 
    \begin{cases}
    1,  \text{ if} \hspace{1em} \tau_i + \Delta_i + f_{ij} \leq \tau_j, \hspace{1em} M_i \cap M_j \neq \emptyset  \\
    0, \text{ otherwise}
    \end{cases}
    & & i,j \in D \\
    B_{sd} &=
    \begin{cases}
    1, \text{ if} \hspace{1em}  f_{sd} \leq \tau_d, s \in S^r,r \in M_d \\
    0, \text{ otherwise}
    \end{cases}
    &  & s \in S^r, r \in R, d \in D
\end{align*}   
Note that since arcs in our directed graph were only created if demands could be initially served by a resource, or served back-to-back by a resource, then any path in the graph is a feasible schedule for a resource that ensures all demands are served on time and for their service duration. 

\begin{definition}
A feasible solution for $mRmD$ is defined as a set of paths from which it is easy to build in polynomial time the set of demands that is met.
\end{definition}

\section{Integer Programming Formulation}\label{formulation_section}
\begin{table}[!h]
    \centering
    \small\setlength\tabcolsep{4.5pt}
    \caption{Notation for problem description and IPM}
    \begin{tabular}{l l}
    \hline 
    \hline
    $D$ & Set of demands \\
    $R$ & Set of resource types \\
    $S^r$ & Set of starting locations for $r \in R$ \\
    $l_s^r$ & Number of resources of type $r \in R$ starting at location $s \in S^r$  \\
    $f_{ij}$ & Travel time between $i \in \{\cup_{r \in R}S^r\} \cup D, j \in D$ \\
    $M_d$ & Subset of resources required by $d \in D$\\
           $\tau_d$ & Service start time for $d \in D$\\
   $\Delta_d$ & Service duration for  $d \in D$\\ $w_d$ & Reward for meeting $d \in D$\\
    \hline
    \hline
    $\boldsymbol{y_d}$= & $\begin{cases}
        1,  \text{ if demand } d \text{ is satisfied} \\
    0, \text{ otherwise} \hspace{12em} d \in D
    \end{cases}$\\
    \\
    $\boldsymbol{x_{ij}^r}$= & $\begin{cases}
        
        1, \text{ if resource type } r \text{ serves demand } \\
       
        \hspace{2em} j \text{ after demand } i \\
        
        0, \text{ otherwise} \hspace{12em} i,j,r: A_{ij}=1, r \in M_i \cap M_j, i,j \in D \\
    \end{cases}$\\
    \\
    $\boldsymbol{x_{sd}^r}$= & $\begin{cases}
        
        1, \text{ if resource type } r \\
            
        \hspace{2em} \text{from starting location } s \\
            
        \hspace{2em} \text{initially serves demand } d \\
            
        0, \text{ otherwise} \hspace{12em} s,d,r: B_{sd}=1, s \in S^r, r \in M_d, d \in D \\
    \end{cases}$\\
    \\
    $\boldsymbol{x_{dt}^r}$= & $\begin{cases}
        
        1, \text{ if resource type } r \\
            
        \hspace{2em} \text{serves demand } d \text{ last} \\
        
        0, \text{ otherwise} \hspace{12em} r \in M_d, d \in D \\
    \end{cases}$\\
    \hline
    \hline
    \end{tabular}
    \label{formulation}
\end{table}

In this section, we present an integer programming model (IPM) for $mRmD$. For an overview of the notation and description of the decision variables, please refer to Table \ref{formulation}.
Given our description of a an instance of $mRmD$ and its feasible solution in Section \ref{prob_desc}, $x$ variables describe the arcs traversed by resources and $y$ denotes whether or not a demand node was visited by all of its required resource types. Note that we use $t$ to denote a dummy sink.
\begin{align}
    & \max ~ \sum_{d \in D}y_dw_d &
\end{align}

s.t.
\begin{align}
      \sum_{d \in D|B_{sd}=1}x_{sd}^r \leq l_s^r &  & s \in S^r, r \in R \\
    \sum_{s \in S^r|B_{si}=1} x_{si}^r+ \sum_{h \in D|A_{hi}=1} x_{hi}^r = \sum\limits_{j \in D|A_{ij}=1} x_{ij}^r + x_{it}^r &  & r \in M_i, i \in D \\
    \sum_{s \in S^r|B_{si}=1} x_{si}^r+ \sum_{h \in D|A_{hi}=1} x_{hi}^r \geq y_i &  & r \in M_i, i \in D \\
    y_d \in \{0,1\} &  & d \in D \\
    x_{sd}^r \in \{0,1\} &  & B_{sd} = 1,s \in S^r, r \in M_d, d \in D \\
    x_{ij}^r \in \{0,1\} &  & A_{ij}=1,r \in M_i \cap M_j, i, j \in D \\
     x_{dt}^r \in \{0,1\} &  & r \in M_d, d \in D
\end{align}

The objective function (1) of IPM is to maximize the overall reward obtained from satisfying demands. Constraints (2) ensure that the number of resources that leave   starting location $s$ is less than or equal to the number of resources available at that location. Constraints (3) maintain flow conservation, i.e.,  the flow of resources arriving to and leaving a demand must be equal. Constraints (4) enforce that a demand is satisfied only if it has been serviced by all of its required resource types. Domain constraints for the variables are given in (5)-(8).

Another valid formulation is to consider time indices and define variables for whether a specific resource of a certain type is at a given location at each time index (see Appendix \ref{alt_form_IPM} for details). However, this formulation could involve significantly more variables and constraints than IPM and is computationally prohibitive, as shown in Table \ref{IPM_vs_Alt-IPM} of Appendix \ref{alt_form_IPM}.
\section{Complexity Results}\label{complexity_results}
First, we show complexity results for $mRmD$, including flow decomposition techniques to determine the feasibility of serving all demands in an $mRmD$ instance. Then we show that $mRmD$ cannot be approximated within a certain order. Finally, we show complexity results for special cases of $mRmD$, which are then used in later sections. Tables \ref{complexity} and \ref{complexity_summary} summarize $mRmD$ and some of its special cases and corresponding complexity results shown in this section.  
\begin{table}[!htb]
\centering
\small\setlength\tabcolsep{3pt}
\caption{Special Cases of the Problem}
\begin{tabular}{|l|l|}
\hline
\textbf{Name}& \textbf{Description} \\ \hline
$1R1D$ & Single resource type; this is equivalent to $1RmD$. \\ \hline
$mR1D$ & Multiple resource types, each demand requires a single resource type. \\ \hline
$2RmD$ & Two resource types, each demand requires either one or both resource types. \\ \hline
$mR\{\text{1 or all}\}D$ & Multiple resource types, each demand requires either one or all resource types. \\ \hline
$mR\{\text{1 or 2}\}D$& Multiple resource types, each demand requires either one or two resource types. \\ \hline
$mRmD$ & \begin{tabular}[c]{@{}l@{}} Multiple resource types, each demand requires a subset. This is the general case.\end{tabular} \\ \hline
\end{tabular}
\label{complexity}
\end{table}
\begin{table}[!htb]
\centering
\small\setlength\tabcolsep{4.5pt}
\caption{Summary of Complexity Results}
\begin{tabular}{|c|c|c|c|c|c|c|}
\hline
\textbf{Name}  & \textbf{\begin{tabular}[c]{@{}c@{}}Travel\\ Time\end{tabular}} & \textbf{\begin{tabular}[c]{@{}c@{}}Service\\ Time\end{tabular}} & \textbf{\begin{tabular}[c]{@{}c@{}}Demand\\ Start Time\end{tabular}} & \textbf{Theorem} & \textbf{Complexity} \\ \hline
$1R1D$  & any   & any  & any   & \ref{1rt-1d} & P \\ \hline
$mR1D$ & any & any  & any& \ref{1rt-1d}  & P \\ \hline
$2RmD$  & 0 & $\infty$ & 0 & \ref{m_1orall} & P \\ \hline
$2RmD$ & 0 & 1 & any & \ref{mrt_1orall_1}  & P \\ \hline
$2RmD$  & any & any & any & \ref{gen_hard} & NP-hard \\ \hline
$mR\{\text{1 or all}\}D$ & 0 & $\infty$ & 0 & \ref{m_1orall} & P \\ \hline
$mR\{\text{1 or all}\}D$ & 0 & 1 & any & \ref{mrt_1orall_1}  & P  \\ \hline
$mR\{\text{1 or 2}\}D$ & any & any & any & \ref{mrt-1or2} & NP-hard \\ \hline
$mRmD$ & any & any & any & \ref{gen_hard},\ref{non_approx} & NP-hard \\ \hline
\end{tabular}
\label{complexity_summary}
\end{table}

\begin{theorem}

\label{gen_hard}
$mRmD$ is NP-hard, for $|R| \geq 2$.
\end{theorem}

\proof{Proof.} 
Refer to Appendix \ref{pf_2rt-md} for the details.
\endproof

\subsection{Flow Decomposition}
\label{flow_decomp_section}
In this subsection, we introduce flow decomposition techniques to determine the feasibility of serving all demands in $mRmD$ instances. We first define some additional notation that is used here, as well as in Section \ref{algorithms}. Let $|R|$ be the number of resource types (labeled $1,2,\ldots,|R|)$ and $D^r \subseteq D$ be the set of demands that require a resource unit of type $r \in R$, i.e., $D^r \coloneqq \{i \in D| r \in M_i\}$.

\begin{theorem}
\label{flow_decomp}
The convex hull induced by the feasible points of IPM remains the same when the $x$ variables are relaxed.
\end{theorem}

\proof{Proof.} Consider the convex hull induced by the feasible points of IPM. We show that when the $x$ variables are relaxed, the extreme points remain the same. Let the binary $y$ variables be fixed (i.e., the demands met are known). In the IPM, when $y$ variables are fixed, the objective function (1) becomes fixed, constraints (5) are removed, and the formulation can be decomposed into $|R|$ maximization problems as follows
\begin{align*}
    \max ~& 0 &
    \tag{1a}
\end{align*}

s.t.
\begin{align*}
      \sum_{d \in D^r|B_{sd}=1}x_{sd}^r \leq l_s^r &  & s \in S^r \tag{2a} \\
    \sum_{s \in S^r|B_{si}=1} x_{si}^r+ \sum_{h \in D^r|A_{hi}=1} x_{hi}^r = \sum\limits_{j \in D^r|A_{ij}=1} x_{ij}^r + x_{it}^r & & i \in D^r \tag{3a}\\
    \sum_{s \in S^r|B_{si}=1} x_{si}^r+ \sum_{h \in D^r|A_{hi}=1} x_{hi}^r \geq y_i & & i \in D^r \tag{4a}\\
    x_{sd}^r \in \{0,1\} &  & B_{sd}=1,s \in S^r, d \in D^r \tag{6a}\\
    x_{ij}^r \in \{0,1\} &  &A_{ij}=1,i, j \in D^r \tag{7a}\\
     x_{dt}^r \in \{0,1\} &  & d \in D^r\tag{8a}
\end{align*}

We denote this formulation as IPM-R; note that it is a feasibility problem.

\begin{lemma}
\label{IPMR_MCF}
IPM-R can be modeled and solved as a MCF.
\end{lemma}

\proof{Proof}
See Appendix \ref{lemma_ipmr} for the details.
\endproof

Since IPM-R can be modeled and solved as a MCF (Lemma \ref{IPMR_MCF}), then the $x$ variables can be relaxed and maintain the same extreme points.
\Halmos
\endproof
Note that when all demands can be satisfied in an instance of $mRmD$, this is equivalent to setting all the $y$ variables to 1 and modeling the instance as $|R|$ feasibility problems that can be solved in polynomial time.
\begin{corollary}
\label{flow_decomp_cor}
The feasibility of serving all demands of an $mRmD$ instance can be determined in polynomial time.
\end{corollary}

Next we show that it is hard to even approximately solve $mRmD$ under natural complexity theoretic assumptions.
\begin{theorem}
\label{non_approx}
$mRmD$ cannot be approximated within $O(\min\{|R|^{\frac{1}{2}-\epsilon},|D|^{1-\epsilon}\})$ unless $NP \subseteq ZPP$ (for every $\epsilon > 0$), where $|R|$ is the number of resource types and $|D|$ is the number of demands.
\end{theorem}

\proof{Proof.} See Appendix \ref{pf_non_approx}.
\endproof

For the remainder of this section, we present complexity results for special cases of $mRmD$ showing that not only the general problem but even very special cases remain NP-hard. Moreover, we identify other special cases that are solvable in polynomial time.
\begin{theorem}
\label{mrt-1or2}
$mR\{\text{1 or 2}\}D$ is NP-hard. 
\end{theorem}

\proof{Proof.} Refer to Appendix \ref{pf_mrt-1or2} for the details.
\endproof

\begin{theorem}
\label{m_1orall}
The special case of $mR\{\text{1 or all}\}D$ where travel times are zero, service times are infinite, and all demands have the same service start time is solvable in polynomial time.
\end{theorem}

\proof{Proof.} See Appendix \ref{pf_m_1orall}.
\endproof

\begin{theorem}

\label{mrt_1orall_1}
The special case of $mR\{\text{1 or all}\}D$ where travel times are zero and service times are one is solvable in polynomial time.
\end{theorem}

\proof{Proof.} Since all service times are one and service start times are integral, then we can group demands by their start times and apply Algorithm \alge $ $ (introduced in Appendix \ref{pf_m_1orall}) to each independent group of demands.
\Halmos
\endproof

\begin{theorem}
\label{1rt-1d}
$1R1D$ and $mR1D$ are solvable in polynomial time.
\end{theorem}

\proof{Proof.}$1R1D$ can be modeled as a minimum cost flow (MCF) problem on a directed acyclic graph. See Appendix \ref{pf_1rt-1d}.  $mR1D$ can be decomposed into $m$ independent $1R1D$ problems, which are solvable in polynomial time. \Halmos
\endproof

MCF problems can be  solved in polynomial time (\cite{mcf_ref}). For example, the minimum mean cycle-canceling algorithm is a strongly polynomial algorithm which runs in $O(m^2n\log n)$ time, where $n,m$ are the number of nodes and arcs, respectively, in the network (\cite{cycle_canceling}).

\textit{Observation 1.} If the travel times are zero, then $1R1D$ can be formulated as a maximum weighted coloring problem on an interval graph. See Appendix \ref{obs_1}.

\textit{Observation 2.} If the travel times are not zero, $1R1D$ cannot be formulated as a maximum weighted coloring problem on an interval graph. This can been shown simply through the use of counterexamples and the characterizations of interval graphs (\cite{gilmore}).

\section{Approximation Algorithms}\label{algorithms}
In this section, we present approximation algorithms for special cases of $mRmD$ and prove performance guarantees. Table \ref{approx_summary_table} presents a summary of the results.

\begin{table}[!htb]
\centering
\caption{Summary of Approximation Algorithm Results}
\begin{tabular}{|c|c|c|c|}
\hline
\textbf{Name} & \multicolumn{1}{l|}{\textbf{Theorem}} & \textbf{Result} & \textbf{\begin{tabular}[c]{@{}c@{}}Resource\\ Starting Location\end{tabular}} \\ \hline
Algorithm \alga & \ref{r_approx} & $|R|$-approx.  & same \\ \hline
Algorithm \algb & \ref{chi_approx} & $\chi$-approx.  & same \\ \hline
Algorithm \algc & \ref{ab_approx} & $\frac{a}{b}$-approx.  & same \\ \hline
Bicriteria & \ref{bicriteria} & $\Big((1-\frac{1}{k}),(\frac{1}{1-k\epsilon})\Big)$-bicriteria approx.  & any \\ \hline
\end{tabular}
\label{approx_summary_table}
\end{table}
In the remainder of this section, we use the following notation.
Let $m_i$ be the number of resources of type $i$ for $i = 1,\ldots,|R|$. Without loss of generality, we assume that $m_1 \leq m_2 \leq \ldots \leq m_{|R|}$. Let sub-problem $1R1D^r$ refer to a $1R1D$ problem instance that arises from the $mRmD$ instance at hand by considering only the $r^{th}$ resource type ($m_r$ resources) and the corresponding subset of relevant demands $D^r$, as defined in Section \ref{flow_decomp_section}. 
Let $OPT(\mathcal{I})$ represent the optimal objective value for instance $\mathcal{I}$ of $mRmD$. Similarly, let $\chi(\mathcal{I})$ represent the objective value of the solution produced by applying Algorithm $\chi$ to instance $\mathcal{I}$ of $mRmD$.

\begin{definition}
An $\alpha$-approximation algorithm for $mRmD$, where $\alpha > 1$, is an algorithm that for every instance of $mRmD$, returns a feasible solution with objective value at least $\frac{1}{\alpha}$ times the optimal objective. Moreover the running time of the algorithm is polynomial in the size of the instance.
\end{definition}

\begin{algorithm}[!htb]
\renewcommand{\thealgocf}{}
\caption{\alga}
\label{r_approximation}
 \SetKwInOut{Input}{Input}
 \SetKwInput{Output}{Output}
 \Input{$\mathcal{I}$}
 \alga$(\mathcal{I}) \leftarrow 0$ \\
  \For{$1 \leq r \leq |R|$}{
    Solve $1R1D^r$; Let \alga$_r(\mathcal{I})$ denote the objective value.  \\
    \If{\alga$_r(\mathcal{I}) \geq$ \alga$(\mathcal{I})$}{
        \alga$(\mathcal{I}) \leftarrow$\alga$_r(\mathcal{I})$ \\
    }
    $D \leftarrow D \setminus D^r$
  }
  \Output{\alga$(\mathcal{I})$}
\end{algorithm}
\textbf{Algorithm \alga $ $ Description:} Algorithm \alga $ $ iterates through the $|R|$ resource types and solves $1R1D^r$ to determine the value of each resource type. The algorithm then picks the best resource type found and serves the optimal set of demands that request that type.

\begin{theorem}
\label{r_approx}
Algorithm \alga $ $ is an $|R|$-approximation algorithm for $mRmD$ instances with the same starting location for every resource.
\end{theorem}

\proof{Proof.} 
First note that using Algorithm \alga $ $, we may construct a feasible solution to instance $\mathcal{I}$ of $mRmD$ with an objective value of \alga($\mathcal{I})$. Because all resources have the same starting location and demands $D^r$ only request resource type $r$ or larger (where larger resource types have at least as many resources as type $r$), then all other resource types larger than $r$ may follow the same schedule determined by the solution to $1R1D^r$ to ensure all demands are met by their requested resource types. Thus, the feasible solution will consist of the paths from the $1R1D^r$ sub-problem that produced the maximum objective value for all resource types greater than or equal to $r$.

Now, we may partition the optimal objective value $OPT(\mathcal{I})$ for instance $\mathcal{I}$ of $mRmD$ into values $OPT_r(\mathcal{I}), r= 1,\ldots,|R|$, such that $\sum_{i=1}^{|R|}OPT_r(\mathcal{I})=OPT(\mathcal{I})$ and each $OPT_r(\mathcal{I})$ represents the objective value of all demands met whose smallest resource type requests is $r$.

\begin{lemma}
\alga$_r(\mathcal{I}) \geq OPT_r(\mathcal{I})$ for $r = 1,\ldots,|R|$.
\end{lemma}

\proof{Proof.}
Clearly, $OPT_r(\mathcal{I})$ is a feasible objective value for the $r^{th}$ iteration of Algorithm \alga $ $, as at the end of each iteration, we remove $D^r$ from $D$. That is, the sets of demands considered in each iteration of Algorithm \alga $ $ are also partitioned according to the indices of their requested resource types. However, \alga$_r(\mathcal{I})$ is the optimal objective value of the $r^{th}$ iteration of Algorithm \alga $ $ and so \alga$_r(\mathcal{I}) \geq OPT_r(\mathcal{I})$.\Halmos
\endproof

Further, \alga$(\mathcal{I}) \geq $\alga$_r(\mathcal{I})$ for $r = 1, \ldots, |R|$ and so we have
$$\sum_{r=1}^{|R|}\text{\alga}(\mathcal{I}) \geq \sum_{r=1}^{|R|}\text{\alga}_r(\mathcal{I}) \geq \sum_{r=1}^{|R|}OPT_r(\mathcal{I}) = OPT(\mathcal{I}).$$
This implies
$$|R|\cdot \text{\alga}(\mathcal{I}) \geq OPT(\mathcal{I})$$
and so, dividing both sides by $|R|$, we obtain $$\text{\alga}(\mathcal{I}) \geq \frac{OPT(\mathcal{I})}{|R|}. \Halmos$$ 
\endproof

\begin{remark}
Note that Algorithm \alga $ $ can be solved in polynomial time as it relies on solving $|R|$ MCF problems.
\end{remark}

For a problem instance $\mathcal{I}$ of $mRmD$, we create a conflict graph $G = (V,E)$ as follows: \\
\textbf{Nodes} \\
Create a node $v=M_j$ for each unique $M_j$, $j \in D$ \\
\textbf{Arcs} \\
$(u,v) \in E$ if $u \cap v \neq \emptyset$, $u,v \in V$  \\
Let $\bar{D}^v \coloneqq \{d \in D| M_d = v\}$ for all $v \in V$ and the sub-problem $1R1\bar{D}^v$ refer to a $1R1D$ problem instance that arises from the $mRmD$ instance at hand by considering all resources of type $i$ where $i$ is the smallest index in set $v$ and the corresponding subset of relevant demands $\bar{D}^v$.

\begin{algorithm}[!htb]
\renewcommand{\thealgocf}{}
\caption{\algb}
\label{chi-approximation}
 \SetKwInOut{Input}{Input}
 \SetKwInput{Output}{Output}
 Suppose there exists a $\chi$-coloring of $G$. Label the colors $1,2,\ldots,\chi$. Let $S_i$ be the set of nodes that are colored $i$. \\
 \Input{$\mathcal{I},G,\chi$, $S_i$, $i \in \{1,2,\ldots,\chi\}$}
 \algb$(\mathcal{I}) \leftarrow 0$\\
  \For{$1 \leq i \leq \chi$}{
    \For{$v \in S_i$}{
        Solve $1R1\bar{D}^v$; Let \algb$_v(\mathcal{I})$ denote the objective value.\\
    }
    \If{$\sum_{v\in S_i} \text{\algb}_v(\mathcal{I}) \geq \text{\algb}(\mathcal{I}):$}{
        $\text{\algb}(\mathcal{I}) \leftarrow \sum_{v\in S_i} \text{\algb}_v(\mathcal{I})$
    }
  }
  \Output{$\text{\algb}(\mathcal{I})$}
\end{algorithm}

\textbf{Algorithm \algb $ $ Description:} Given a $\chi$-coloring of the conflict graph $G$, Algorithm \algb $ $ iterates through the $\chi$ colors in the coloring of $G$ and solves $1R1\bar{D}^v$ for all nodes $v$ of that color to determine the value of each color. The algorithm then picks the best color found and serves the optimal set of demands of that color.

\begin{theorem}
\label{chi_approx}
Given a $\chi$-coloring of the conflict graph $G$, Algorithm \algb $ $ is a $\chi$-approximation algorithm for $mRmD$ instances with the same starting location for every resource.
\end{theorem}

\proof{Proof.}
First note that using Algorithm \algb $ $, we may construct a feasible solution to instance $\mathcal{I}$ of $mRmD$ with an objective value of \algb($\mathcal{I}$). Because all resources have the same starting location and the resource type used to solve $1R1\bar{D}^v$ is the one with the smallest number of resources in the subset $v$, then all other resource types in $v$ may follow the same schedule determined by the solution to $1R1\bar{D}^v$ to ensure all demands are met by their requested resource types. Thus, the feasible solution will consist of the paths from the $1R1\bar{D}^v$ sub-problem for all resources in $v$ and for all $v \in S_i$, for the color $i$ that produced the maximum objective value.

Now, we may partition the optimal objective value $OPT(\mathcal{I})$ for instance $\mathcal{I}$ of $mRmD$ into values $OPT_v(\mathcal{I}), v \in S_i, i \in \{1,\ldots,\chi\}$, where $\sum_{i=1}^{\chi}\sum_{v\in S_i}OPT_v(\mathcal{I}) = OPT(\mathcal{I})$ and each $OPT_v(\mathcal{I})$ represents the objective value of all demands $d$ met whose subset request $M_d = v$.

\begin{lemma}
$\text{\algb}_v(\mathcal{I}) \geq OPT_v(\mathcal{I})$, $v \in S_i$, $i = 1,\ldots,\chi$.
\end{lemma}

\proof{Proof.} Clearly, $OPT_v(\mathcal{I})$ is a feasible objective value to $1R1\bar{D}^v$. However, $\text{\algb}_v(\mathcal{I})$ is the optimal objective value to $1R1\bar{D}^v$ and so, $\text{\algb}_v(\mathcal{I}) \geq OPT_v(\mathcal{I})$. \Halmos
\endproof

Further, $\text{\algb}(\mathcal{I}) \geq \sum_{v \in S_i}\text{\algb}_v(\mathcal{I})$ for $i = 1, \ldots, \chi$ and so we have
$$\sum_{i=1}^{\chi}\text{\algb}(\mathcal{I}) \geq \sum_{i=1}^{\chi}\sum_{v\in S_i}\text{\algb}_v(\mathcal{I}) \geq \sum_{i=1}^{\chi}\sum_{v\in S_i}OPT_v(\mathcal{I}) = OPT(\mathcal{I}).$$
This implies
$$\chi \cdot \text{\algb}(\mathcal{I}) \geq OPT(\mathcal{I})$$
and so, dividing both sides by $\chi$, we obtain
$$\text{\algb}(\mathcal{I}) \geq \frac{OPT(\mathcal{I})}{\chi}. \Halmos$$
\endproof

\begin{remark}
Note that given a $\chi$-coloring of the conflict graph $G=(V,E)$, Algorithm \algb $ $ can be solved in polynomial time as it relies on solving $|V|$ MCF problems.
\end{remark}

Theorem \ref{chi_approx} may be generalized to an $a:b$-coloring. 
\begin{definition}
An $a:b$-coloring is a $b$-fold coloring out of $a$ available colors, where a $b$-fold coloring is an assignment of sets of size $b$ to nodes of a graph such that adjacent nodes receive disjoint sets.
\end{definition}

\begin{algorithm}[!htb]
\renewcommand{\thealgocf}{}
\caption{\algc}
\label{ab-algo}
 \SetKwInOut{Input}{Input}
 \SetKwInput{Output}{Output}
Suppose there exists an a:b-coloring of $G$. Label the colors $1,2,\ldots,a$. Let $S_i$ be the set of nodes that are fractionally colored $i$. \\
 \Input{$\mathcal{I},G,a$, $S_i$, $i \in \{1,2,\ldots,a\}$}
 $\text{\algc}(\mathcal{I}) \leftarrow 0$\\
  \For{$1 \leq i \leq a$}{
    \For{$v \in S_i$}{
        Solve $1R1\bar{D}^v$; Let $\text{\algc}_v(\mathcal{I})$ denote the objective value.
    }
    \If{$\sum_{v\in S_i} \text{\algc}_v(\mathcal{I}) \geq \text{\algc}(\mathcal{I}):$}{
        $\text{\algc}(\mathcal{I}) \leftarrow \sum_{v\in S_i} \text{\algc}_v(\mathcal{I})$
    }
  }
  \Output{$\text{\algc}(\mathcal{I})$}
\end{algorithm}

\textbf{Algorithm \algc $ $ Description:} Given an $a:b$-coloring of the conflict graph $G$, Algorithm \algc $ $ iterates through the $a$ colors in the $a:b$-coloring of $G$ and solves $1R1\bar{D}^v$ for all nodes $v$ of that color to determine the value of each color. The algorithm then picks the best color found and serves the optimal set of demands of that color.

\begin{theorem}
\label{ab_approx}
Given an $a:b$-coloring of the conflict graph $G$, Algorithm \algc $ $ is an $\frac{a}{b}$-approximation algorithm for $mRmD$ instances with the same starting location for every resource.
\end{theorem}

\proof{Proof.}
First note that using Algorithm \algc $ $, we may construct a feasible solution to instance $\mathcal{I}$ of $mRmD$ with an objective value of \algc($\mathcal{I}$). Because all resources have the same starting location and the resource type used to solve $1R1\bar{D}^v$ is the one with the smallest number of resources in the subset $v$, then all other resource types in $v$ may follow the same schedule determined by the solution to $1R1\bar{D}^v$ to ensure all demands are met by their requested resource types. Thus, the feasible solution will consist of the paths from the $1R1\bar{D}^v$ sub-problem for all resources in $v$ and for all $v \in S_i$, for the color $i$ that produced the maximum objective value.

Now, we may partition the optimal objective value $OPT(\mathcal{I})$ for instance $\mathcal{I}$ of $mRmD$ into values $OPT_j(\mathcal{I}), j = 1,\ldots,a$, where $b\cdot OPT(\mathcal{I}) = \sum_{j=1}^{a}OPT_j(\mathcal{I})$ and each $OPT_j(\mathcal{I})$ represents the objective value of all demands met that are fractionally colored $j$.

\begin{lemma}
$\sum_{v\in S_i}\text{\algc}_v(\mathcal{I}) \geq OPT_i(\mathcal{I})$ for all $i = 1,\ldots,a$.
\end{lemma}

\proof{Proof.}
Clearly, $OPT_i(\mathcal{I})$ is a feasible objective value to the $i^{th}$ iteration of Algorithm \algc $ $. However, $\sum_{v\in S_i}\text{\algc}_v(\mathcal{I})$ is the optimal objective value to the $i^{th}$ iteration and so, $\sum_{v\in S_i}\text{\algc}_v(\mathcal{I}) \geq OPT_i(\mathcal{I})$. \Halmos
\endproof

Further, $\text{\algc}(\mathcal{I}) \geq \sum_{v \in S_i}\text{\algc}_v(\mathcal{I})$ for $i = 1, \ldots, a$ and so we have
$$\sum_{i=1}^{a}\text{\algc}(\mathcal{I}) \geq \sum_{i=1}^{a}\sum_{v\in S_i}\text{\algc}_v(\mathcal{I}) \geq \sum_{i=1}^{a}OPT_i(\mathcal{I}) = b \cdot OPT(\mathcal{I}).$$
This implies
$$a\cdot \text{\algc}(\mathcal{I}) \geq b \cdot OPT(\mathcal{I})$$
and so, dividing both sides by $a$, we obtain
$$\text{\algc}(\mathcal{I}) \geq \frac{b}{a}\cdot OPT(\mathcal{I}).\Halmos$$
\endproof

\begin{remark}
Note that given an $a:b$-coloring of the conflict graph $G=(V,E)$, Algorithm \algc $ $ can be solved in polynomial time as it relies on solving $b\cdot|V|$ MCF problems.
\end{remark}

\begin{definition}
An $(\alpha,\beta)$-bicriteria approximation algorithm for $mRmD$, $( \alpha > 1, \beta > 1)$, is an algorithm that given any instance of the problem, returns a solution whose objective value is at least $\frac{1}{\alpha}$ fraction of the optimal objective and uses at most $\beta$ times more resources for every resource type. Moreover, the running time of the algorithm is polynomial in size of the input instance.
\end{definition}
In the following theorem, we show that there exists a good bicriteria algorithm for instances where the optimal solution is able to satisfy nearly all the demand. The result can be interpreted as a smooth degradation of Corollary~\ref{flow_decomp_cor} that shows the decision problem of deciding whether all demands are satisfiable is polynomial time solvable.

\begin{definition}
Given a parameter $\epsilon > 0$, an instance of $mRmD$ is $(1-\epsilon)$-satisfiable if the optimal solution has objective at least $(1-\epsilon)$ times the total reward possible, $\sum_{d\in D} w_d$.
\end{definition}

\begin{theorem}
\label{bicriteria}
For any $0 < \epsilon < 1$ and $1\leq k\leq \frac{1}{\epsilon}$, there exists a $\Big(\frac{k-1}{k},\frac{1}{1-k\epsilon}\Big)$-bicriteria approximation algorithm for $(1-\epsilon)$-satisfiable instances of $mRmD$.
\end{theorem} 

\begin{remark}
Note that when we multiply the number of resources of each resource type by $(\frac{1}{1-k\epsilon})$, there is a possibility that the number of resources of each type is no longer integral. Without loss of generality, we may round up any non-integral values.
\end{remark}

\begin{remark}
As an instantiation of Theorem~\ref{bicriteria}, consider the parameter $\epsilon=0.01$ and $k=5$. Then, given an instance where the optimal solution satisfies at least $99\%$ of the weighted sum of demands, Theorem~\ref{bicriteria} returns a solution that uses 25\% more resources of every resource type whose objective value is at least 94\% of the weighted sum of demands.
\end{remark}

\proof{Proof of Theorem \ref{bicriteria}.}
Clearly, the objective of the linear relaxation, $OPT_{LP}(\mathcal{I})$, of IPM (i.e., all variables are within the range $[0,1]$) for $(1-\epsilon)$-satisfiable instances of $mRmD$ must also be at least $(1-\epsilon)$ times the total reward possible, $\sum_{d\in D}w_d$. The algorithm uses the solution to the linear relaxation of IPM and Algorithm $\mathcal{D}$ to construct an integral solution to IPM whose objective is at least $(\frac{k-1}{k})OPT_{LP}(\mathcal{I})$ and uses no more than $(\frac{1}{1-k\epsilon})$ times more resources for every resource type. Let $(\mathbf{y}^{LP},\mathbf{x}^{LP})$ represent the optimal solution to the linear relaxation and define the following sets:
$$B = \{i: y_i^{LP} < 1- k\epsilon\} \text{ and } G = D \setminus B.$$
Now, consider the solution $(\mathbf{y}',\mathbf{x}')$ where 
\begin{align*}
y_i'&=\begin{cases}
1 \hspace{1em} \forall i \in G\\
0 \hspace{.8em} \text{ else, } \end{cases} 
& \\
x_{ij}^{r'} &= \frac{x_{ij}^{rLP}}{1-k\epsilon} & r \in M_i \cap M_j, i,j \in D, \\
x_{sd}^{r'} &= \frac{x_{sd}^{rLP}}{1-k\epsilon} & s \in S^r, r \in M_d, d \in D, \\
x_{dt}^{r'} &= \frac{x_{dt}^{rLP}}{1-k\epsilon} & r \in M_d, d \in D.
\end{align*}

Note that since $x$ variables may be relaxed (see Theorem \ref{flow_decomp}), $(\mathbf{y}',\mathbf{x}')$ is a solution to IPM when resources of each resource type are increased by a factor of $(\frac{1}{1-k\epsilon})$ that satisfies all constraints except the upper bound on the $x$ variables.
\begin{lemma}
\label{sol_transform}
A solution $(\mathbf{y},\mathbf{x})$ to IPM in which the $x$ variables are unbounded above may be transformed into a feasible solution $(\mathbf{y},\mathbf{\bar{x}})$ to IPM in which the $x$ variables are bounded above by 1. 
\end{lemma} 
\proof{Proof.}
Note that from Section \ref{prob_desc}, we may describe the solution $(\mathbf{y},\mathbf{x})$ as a directed acyclic graph $G = (V,A)$ where each $x^r_{ij} \in \mathbf{x}$ corresponds to the flow along arc $(i,j) \in A$ for resource type $r \in R$ and $i,j \in V$. Now, the unbounded solution $(\mathbf{y},\mathbf{x})$ may be transformed into the bounded solution $(\mathbf{y},\mathbf{\bar{x}})$ by the procedure described in Algorithm \algd $ $. The number of operations in Algorithm \algd $ $ is finite since at each stage, the total flow along all edges in $G$ is decreasing as we send the same flow value on fewer number of arcs. It is also important to note that all variables used in Algorithm \algd $ $ must exist by the triangle inequality. Further, the $\mathbf{y}$ variable values remain valid because for every node in which flow in and flow out is reduced, the flow is never reduced to less than 1. \Halmos
\begin{algorithm}[!htb]
\renewcommand{\thealgocf}{}
\caption{\algd}
\label{algo_sol_trans}
 \SetKwInOut{Input}{Input}
 \SetKwInput{Output}{Output}
 \Input{$(\mathbf{y},\mathbf{x})$}
 $\mathbf{\bar{x}} \leftarrow \mathbf{x}$ \\
\While{There are still $\bar{x}^r_{ij} > 1$}{
\If{$i$ is a source and $j$ is a sink}{
Reduce flow value of $\bar{x}^r_{ij}$ to 1.}
\If{$i$ is a source and $j$ is not a sink}{
Pick some $\bar{x}^r_{jk} > 0$. Let $\delta = \min\{\bar{x}^r_{ij}-1,\bar{x}^r_{jk}\}$. Reduce flow value along $\bar{x}^r_{ij}-\bar{x}^r_{jk}$ path by $\delta$ and increase flow value of $\bar{x}^r_{ik}$ by $\delta$.}
\If{$i$ is not a source and $j$ is a sink}{
Pick some $\bar{x}^r_{hi} > 0$. Let $\delta = \min\{\bar{x}^r_{ij}-1,\bar{x}^r_{hi}\}$. Reduce flow value along $\bar{x}^r_{hi}-\bar{x}^r_{ij}$ path by $\delta$ and increase flow value of $\bar{x}^r_{hj}$ by $\delta$.}
\Else{
Pick some $\bar{x}^r_{hi}, \bar{x}^r_{jk} > 0$. Let $\delta = \min\{\bar{x}^r_{ij}-1,\bar{x}^r_{hi},\bar{x}^r_{jk}\}$. Reduce flow value along $\bar{x}^r_{hi}-\bar{x}^r_{ij}-\bar{x}^r_{jk}$ path by $\delta$ and increase flow value of $\bar{x}^r_{hk}$ by $\delta$.}}
\Output{$(\mathbf{y},\mathbf{\bar{x}})$}
\end{algorithm}
\endproof
By Lemma \ref{sol_transform}, we may transform $(\mathbf{y}',\mathbf{x}')$ into a feasible solution to IPM when resources are increased by $(\frac{1}{1-k\epsilon})$ and the $x$ variables satisfy their bound restrictions. Denote this solution by $(\mathbf{y}',\mathbf{x}'')$. 

Now from Theorem \ref{flow_decomp} we know that if the $\mathbf{y}'$ variables are fixed, then IPM decomposes into $|R|$ MCF problems which can be solved in polynomial time and produce an integral optimal solution. Therefore, we can fix the $\mathbf{y}'$ variables and solve IPM optimally, letting $OPT(\mathcal{I}),(\mathbf{y}',\mathbf{x}^*)$ represent the optimal objective value and solution, respectively.

\begin{lemma}
$OPT(\mathcal{I}) \geq (\frac{k-1}{k})OPT_{LP}(\mathcal{I}).$
\end{lemma}

\proof{Proof.}
First note that 
$$\sum_{i \in B}y_i^{LP}w_i + \sum_{i \in G}w_i \geq \sum_{i \in B}y_i^{LP}w_i + \sum_{i \in G}y_i^{LP}w_i = OPT_{LP}(\mathcal{I}) \geq (1-\epsilon)\sum_{i \in B}w_i + (1-\epsilon)\sum_{i \in G}w_i.$$
This implies that 
\begin{align}
    \epsilon\sum_{i \in G}w_i &\geq (1-\epsilon)\sum_{i \in B}w_i - \sum_{i \in B}y_i^{LP}w_i \\
    &> (1-\epsilon)\sum_{i \in B}w_i - (1-k\epsilon)\sum_{i \in B}w_i \\
    &= \epsilon(k-1)\sum_{i \in B}w_i
\end{align}
where (10) comes from the fact that $\sum_{i \in B}y_i^{LP}w_i < (1-k\epsilon)\sum_{i \in B}w_i$. Now, we have that $\sum_{i \in G}w_i \geq (k-1)\sum_{i \in B}w_i$ and so
$$OPT(\mathcal{I}) = \sum_{i \in G}w_i \geq \left(\frac{k-1}{k}\right)\left(\sum_{i \in G}w_i + \sum_{i \in B}w_i\right) \geq \left(\frac{k-1}{k}\right)OPT_{LP}(\mathcal{I}). \Halmos$$
\endproof
Thus, we have shown that $(\mathbf{y}',\mathbf{x}^*)$ is an integral solution to IPM using  $(\frac{1}{1-k\epsilon})$ times more resources for every resource type whose objective value, $OPT(\mathcal{I})$, is greater than or equal to $(\frac{k-1}{k})OPT_{LP}(\mathcal{I})$. \Halmos
\endproof

\section{Extension}\label{prob_ext}
We present an  extension to $mRmD$ (denoted by $mRmDc$) that considers integer-valued travel costs $(c)$ between different locations (e.g., $c_{ij}$ denotes the travel cost between the locations of demands $i$ and $j$ and $c_{sd}$ is the travel cost between a resource's starting location $s$ and the location of demand $d$. The integer programming model incorporating travel costs (IPM-C) is as follows:
\begin{align}
    \max ~ \sum_{d \in D}y_dw_d - \sum_{r \in R}\Big(\sum_{s \in S^r}\sum_{d \in D}x_{sd}^rc_{sd} + \sum_{i \in D}\sum_{j \in D}x_{ij}^rc_{ij}\Big) &
\end{align}

s.t. 
$$(2)-(8)$$

\begin{theorem}
Assume that the travel cost between any two locations is less than or equal to $\frac{1}{2|R|}$ times the minimum reward for a demand $(w^{min})$. Then Algorithm \alga $ $ is a $2|R|$-approximation algorithm for $mRmDc$ instances with the same starting location for every resource.
\end{theorem}

\proof{Proof.}
Note that we may run Algorithm \alga $ $ on instance $\mathcal{I}$ of $mRmDc$ assuming travel costs are zero, where $\text{\alga}(\mathcal{I}),D_{\text{\alga}} \subseteq D$ are the objective value and set of demands met in the solution produced by Algorithm \alga $ $, repsectively. Let $C,C^*$ be the total travel costs incurred from the solutions of Algorithm \alga $ $ and the optimal solution, respectively. Thus, the objective value produced by Algorithm \alga $ $ is $\text{\alga}(\mathcal{I}) - C$ and the optimal objective value is $OPT(\mathcal{I}) - C^*$. Then we have
\begin{align}
    \text{\alga}(\mathcal{I}) - C & \geq \text{\alga}(\mathcal{I}) - c|D_{\text{\alga}}||R| \\
    & \geq \text{\alga}(\mathcal{I}) - \Big(\frac{1}{2|R|}w^{min}\Big)|D_{\text{\alga}}||R| \\
    &\geq \frac{1}{2}\text{\alga}(\mathcal{I}) \\
    &\geq \frac{1}{2|R|}OPT(\mathcal{I})\\
    & \geq \frac{1}{2|R|}(OPT(\mathcal{I}) - C^*).
\end{align}
Note that (14) follows from the assumption about travel costs and (16) is from the results of Theorem \ref{r_approx}. \Halmos
\endproof

\section{Computational Study}\label{comp_study}
To test IPM, we created different problem instances using 2-7 resource types, 100-800 demands, and varied number of resources. In all cases, demand start times are randomly chosen in the interval of [0, 1440], used to represent scheduling a day (in minutes), with demands and resources randomly placed on a 20x20 grid structure. Adjacent nodes on the grid are 1 minute apart. Service times are determined from a triangular distribution with a minimum of 15 minutes, maximum of 120 minutes, and mode of 30 minutes. All resource types have the same number of resources (e.g., if there are 3 resource types and a total of 12 resources then there are 4 resources of each type for that instance). Rewards for demands are proportional to their service time and the number of resource types the demand requires.  Resource requirement subsets are determined randomly, with a 50\% chance that a certain resource type will be required by a demand incident. Results presented are averages of 10 instances. All instances were ran using Gurobi version 8.0.1. 

Tables \ref{results_1} and \ref{results_2} present the run times (in seconds) for all instances considered, where $|R|$ is the number of resource types, $|D|$ is the number of demands, and $L$ is the total number of resources. Note that \textit{Scaled Demands} refers to results for instances in which the reward for demand was multiplied by a factor of 100.  As expected, as the number of demands or number of resource types increases, the run time increases. Computational experiments show that for small and medium sized problems (less than 5 resource types), IPM can be efficiently solved by Gurobi. However, larger sized problems are harder to solve. For example, problem instances with 7 resource types and between 700 and 800 demands took, on average, between 12 and 50 minutes to solve. Results are similar between instances with and without scaling. Appendix \ref{alt_form_IPM} describes an alternate formulation for $mRmD$ and shows that the runtimes for IPM are faster. Appendix \ref{bicriteria_ext} presents a discussion on the trade-off between resource capacity and the objective function value (i.e., demands met), as highlighted by Theorem \ref{bicriteria}.

\begin{table}[!htb]
    \begin{minipage}{0.5\textwidth}
        \centering
        \caption{Results for number of resource types equal to 2,3 and 4}
        \label{results_1}
        \resizebox{!}{.2\paperheight}{%
\begin{tabular}{|c|c|c|c|c|}
\hline
\textbf{$|R|$} & \textbf{$|D|$} & \textbf{L} & \textbf{Run time} & \textbf{\begin{tabular}[c]{@{}c@{}}Run time\\ (Scaled Demands)\end{tabular}} \\ \hline
2 & 100 & 10 & 0.08 & 0.10 \\ \hline
2 & 200 & 14 & 0.25 & 0.33 \\ \hline
2 & 300 & 18 & 0.65 & 0.90 \\ \hline
2 & 400 & 22 & 1.24 & 1.85 \\ \hline
2 & 500 & 26 & 3.10 & 2.74 \\ \hline
2 & 600 & 30 & 4.88 & 5.88 \\ \hline
2 & 700 & 34 & 8.58 & 7.20 \\ \hline
2 & 800 & 38 & 12.06 & 13.96 \\ \hline
3 & 100 & 12 & 0.12 & 0.10 \\ \hline
3 & 200 & 18 & 0.32 & 0.59 \\ \hline
3 & 300 & 24 & 1.44 & 1.48 \\ \hline
3 & 400 & 30 & 2.66 & 3.70 \\ \hline
3 & 500 & 36 & 4.80 & 7.35 \\ \hline
3 & 600 & 42 & 9.04 & 15.42 \\ \hline
3 & 700 & 48 & 22.21 & 21.33 \\ \hline
3 & 800 & 54 & 26.60 & 40.69 \\ \hline
4 & 100 & 16 & 0.15 & 0.25 \\ \hline
4 & 200 & 24 & 0.77 & 1.13 \\ \hline
4 & 300 & 32 & 2.75 & 3.04 \\ \hline
4 & 400 & 40 & 5.71 & 7.89 \\ \hline
4 & 500 & 48 & 18.72 & 16.69 \\ \hline
4 & 600 & 56 & 21.07 & 30.45 \\ \hline
4 & 700 & 64 & 98.42 & 52.30 \\ \hline
4 & 800 & 72 & 147.81 & 72.01 \\ \hline
\end{tabular}}
    \end{minipage}
    \hfill
    \begin{minipage}{0.5\textwidth}
        \centering
        \caption{Results for number of resource types equal to 5,6 and 7}
        \label{results_2}
        \resizebox{!}{.2\paperheight}{%
       \begin{tabular}{|c|c|c|c|c|}
\hline
\textbf{$|R|$} & \textbf{$|D|$} & \textbf{L} & \textbf{Run time} & \textbf{\begin{tabular}[c]{@{}c@{}}Run time\\ (Scaled Demands)\end{tabular}} \\ \hline
5 & 100 & 20 & 0.13 & 0.22 \\ \hline
5 & 200 & 30 & 0.92 & 1.60 \\ \hline
5 & 300 & 40 & 5.39 & 6.31 \\ \hline
5 & 400 & 50 & 8.08 & 17.57 \\ \hline
5 & 500 & 60 & 43.06 & 44.37 \\ \hline
5 & 600 & 70 & 55.82 & 161.77 \\ \hline
5 & 700 & 80 & 145.10 & 100.11 \\ \hline
5 & 800 & 90 & 276.00 & 349.65 \\ \hline
6 & 100 & 24 & 0.24 & 0.40 \\ \hline
6 & 200 & 36 & 1.31 & 2.41 \\ \hline
6 & 300 & 48 & 7.13 & 9.01 \\ \hline
6 & 400 & 60 & 21.74 & 22.78 \\ \hline
6 & 500 & 72 & 87.24 & 76.92 \\ \hline
6 & 600 & 84 & 184.65 & 349.21 \\ \hline
6 & 700 & 96 & 474.84 & 313.35 \\ \hline
6 & 800 & 108 & 1133.85 & 778.51 \\ \hline
7 & 100 & 28 & 0.34 & 0.58 \\ \hline
7 & 200 & 42 & 2.89 & 3.51 \\ \hline
7 & 300 & 56 & 13.07 & 20.42 \\ \hline
7 & 400 & 70 & 111.12 & 65.29 \\ \hline
7 & 500 & 84 & 169.17 & 195.20 \\ \hline
7 & 600 & 98 & 564.02 & 463.19 \\ \hline
7 & 700 & 112 & 754.89 & 985.09 \\ \hline
7 & 800 & 126 & 2965.94 & 2360.52 \\ \hline
\end{tabular}}
    \end{minipage}
\end{table}

\subsection{Computational Study with Travel Costs}
We ran the same instances as above using IPM-C, where travel costs were equivalent to travel times (i.e., demands that are 3 minutes apart have a travel cost of 3 in the objective). Results and analysis can be found in Appendix \ref{comp_ext}.

\section{Conclusions}\label{conclusion}
In this study, we formulated the heterogeneous multi-resource allocation problem ($mRmD$) where each demand requests a subset of resources simultaneously at a specified time, location, and duration as an integer program (IPM). Complexity results were given for $mRmD$, as well as various special cases. A polyhedral result was introduced that allowed us to relax variables in IPM. Further, we developed approximation algorithms for variations of $mRmD$ and proved the correctness of their performance guarantees.  Finally, we tested the performance of the model computationally using Gurobi. 

One simple extension of $mRmD$ (labeled Extension 2) is to consider that resource types may need to visit a destination location before moving on to meet the next demand. For example, after meeting a demand incident, an ambulance may need to drop the patient off at the hospital before proceeding to the next demand location. The problem description and formulation for Extension 2 are provided in Appendix \ref{extension}. As the structure of the solution space does not change significantly, the complexity results in Section \ref{complexity_results} can be applied to this extension, with computational results being very similar to those discussed in Section \ref{comp_study}.  

In this work, we have assumed that all demand requests are known ahead of time (e.g., deterministic). While this formulation can be used to influence planning decisions, future research could consider stochastic and dynamic versions of $mRmD$ (i.e., when all demand incidents may not be known ahead of time). Further, our problem could be generalized to consider requiring multiple units of each resource type.

\ACKNOWLEDGMENT{%
This research has been supported in part by National Science  Foundation (NSF) Graduate Research Fellowship DGE-1650044,  NSF grant CMMI-1538860, NSF- AF:1910423 and NSF-AF:1717947 and the following Georgia Tech benefactors: William W. George, Andrea Laliberte, Joseph C. Mello, Richard “Rick” E. \& Charlene Zalesky, and Claudia \& Paul Raines. The authors would also like to thank the editor and reviewers for their comments and suggestions; their diligent and detailed reviews thoroughly improved our manuscript.
}

%
%
%
\begin{APPENDICES}\label{appendices}
\section{Alternate Formulation for $mRmD$}
\label{alt_form_IPM}
\begin{table}[!h]
    \centering
    \small\setlength\tabcolsep{4.5pt}
    \caption{Notation for IPM-Alt}
    {\begin{tabular}{l l}
    \hline 
    \hline
    $D$ & Set of demands \\
    $R$ & Set of resource types \\
    $S^r$ & Set of starting locations for $r \in R$ \\
    $l_s^r$ & Number of resources of type $r \in R$ starting at location $s \in S^r$  \\
    $\mathcal{K}^r$ & Number of resources of type $r \in R$\\
    $f_{ij}$ & Travel time between $i,j \in \{\cup_{r \in R}S^r\} \cup D$ \\
    $M_d$ & Subset of resources required by $d \in D$\\
           $\tau_d$ & Service start time for $d \in D$\\
   $\Delta_d$ & Service duration for  $d \in D$\\ $w_d$ & Reward for meeting $d \in D$\\
   $T$ & Last time period for demand to be served\\
    \hline
    \hline
    $\boldsymbol{y_d}$= & $\begin{cases}
        1,  \text{ if demand } d \text{ is satisfied} \\
    0, \text{ otherwise} \hspace{12em} d \in D
    \end{cases}$\\
    \\
    $\boldsymbol{x_{it}^{k,r}}$= & $\begin{cases}
        1, \text{ if resource } k \text{ of type } r \text{ is at } \\
        \hspace{2em} \text{ location } i \text{ at time } t \\
        0, \text{ otherwise} \hspace{12em} i,t,k,r: i \in \{\cup_{r \in R}S^r\} \cup D,t \in \{0,\ldots,T\},\\
        \hspace{21em} k \in \{1,\ldots,\mathcal{K}^r\}, r \in R\\
    \end{cases}$\\
       $\boldsymbol{z_{d}^{k,r}}$= & $\begin{cases}
        1, \text{ if demand } d \text{ is met by resource } \\
        \hspace{2em} k \text{ of type }
         r \\
        0, \text{ otherwise} \hspace{12em} d,k,r: k \in \{1,\ldots,\mathcal{K}^r\}, r \in M_d, d \in D\\
    \end{cases}$\\
    \\
    \hline
    \hline
    \end{tabular}}
    \label{alt_formulation}
\end{table}
We present an alternative integer programming model (IPM-Alt) for $mRmD$. For an overview of the notation and description of the decision variables, please refer to Table \ref{alt_formulation}.
{\begin{align}
    & \max ~ \sum_{d \in D}y_dw_d &
\end{align}}

s.t.
{\begin{align}
      \sum_{k=1}^{\mathcal{K}^r}x_{s0}^{k,r} = l_s^r &  & s \in S^r, r \in R \\
    \sum_{s \in S^r}x_{s0}^{k,r} = 1 &  & k \in \{1, \ldots, \mathcal{K}^r\}, r \in R \\
    \sum_{i \in \{\cup_{r \in R}S^r\} \cup D}x_{it}^{k,r} \leq 1 &  & t \in \{0,\ldots,T\},k \in \{1, \ldots, \mathcal{K}^r\}, r \in R \\
    x_{it}^{k,r} + x_{jt'}^{k,r} \leq 1  &  & t' > t |f_{ij} > t' - t, i,j \in \{\cup_{r \in R}S^r\} \cup D, \\
    & & t \in \{0,\ldots,T-1\},k \in \{1,\ldots, \mathcal{K}^r\},r \in R \nonumber\\
    x_{dt}^{k,r} \geq z_{d}^{k,r} & & t \in \{\tau_d,\ldots,\tau_d+\Delta_d\}, k \in \{1, \ldots, \mathcal{K}^r\}, r \in M_d, d \in D \\
    \sum_{k=1}^{\mathcal{K}^r} z_{d}^{k,r} \geq y_d &  & r \in M_d, d \in D \\
    y_d \in \{0,1\} &  & d \in D \\
    x_{it}^{k,r} \in \{0,1\}  &  & i \in \{\cup_{r \in R}S^r\} \cup D,t \in \{0,\ldots,T\},k \in \{1,\ldots,\mathcal{K}^r\}, r \in R\\
    z_{d}^{k,r} \in \{0,1\}  &  & k \in \{1,\ldots,\mathcal{K}^r\}, r \in M_d, d \in D
\end{align}}

The objective function (18) of IPM-Alt is to maximize the overall reward obtained from satisfying demands. Constraints (19) ensure that the number of resources at starting location $s$ at time 0 is equal to the number of resources available at that location. Constraints (20) assign each resource to a starting location. Constraints (21) ensure that a resource can only be at at most one location at each time index. Constraints (22) maintain that the movement of each resource from location to location is valid for the given time indices. Constraints (23) enforce that a demand is met by a given resource only if that resource is at the demand's location for its entire service duration. Constraints (24) state that a demand can only be satisfied if it has been served by a resource of each of its required resource types. Domain constraints for the variables are given in (25)-(27).

We test a small instance of $mRmD$ with 10 demands, 2 resource types (with 3 resources of each type), and demand start times randomly chosen in the interval [0,100]. All other computational details remain the same as described in Section \ref{comp_study}. Table \ref{IPM_vs_Alt-IPM} presents the solution times for 10 runs of this instance using both IPM and IPM-Alt. As shown, even with this very small instance, IPM outperforms IPM-Alt significantly. When attempting to use IPM-Alt to solve the smallest sized instance type described in Section \ref{comp_study} (e.g., 2 resource types, 100 demands, 10 resources), construction of the model did not finish in under an hour due to a large number of variables and constraints.

\begin{table}[!htb]
\caption{Comparison of run times (in seconds)}
\label{IPM_vs_Alt-IPM}
{\begin{tabular}{|c|c|c|c|c|c|c|c|c|c|c|}
\hline
\textbf{Instance} & 1 & 2 & 3 & 4 & 5 & 6 & 7 & 8 & 9 & 10 \\ \hline
\textbf{IPM} & 0.044 & 0.04 & 0.083 & 0.027 & 0.045 & 0.032 & 0.044 & 0.036 & 0.045 & 0.051 \\ \hline
\textbf{IPM-Alt} & 4.47 & 4.15 & 6.33 & 3.57 & 3.58 & 5.82 & 7.22 & 4.20 & 6.19 & 3.19 \\ \hline
\end{tabular}}
\end{table}

\section{Proofs of Theorems and Observations}
\subsection{Proof of Theorem \ref{gen_hard}}
\label{pf_2rt-md}
We prove that $2RmD$ is NP-hard by reduction using Numerical 3-Dimensional Matching (N3DM). The construction is similar to the proof in \cite{aircraft}.\\
\textbf{Instance of N3DM} \\
Integers $t,d \text{ and } a_t,b_t,c_t$ for $i = 1, \ldots, t$, satisfying the following relations: \\
$$\sum_{i=1}^{t}(a_i + b_i + c_i) = td \text{ and } 0 < a_i,b_i,c_i < d \text{ for } i = 1, \ldots, t$$
The goal is to find permutations $\rho$ and $\sigma$ of $\{1, \ldots, t\}$, such that: \\
$$a_i + b_{\rho(i)} + c_{\sigma(i)} = d \text{ for } i = 1, \ldots, t$$

N3DM is shown to be NP-Complete in \cite{garey}. Given an instance of N3DM, we define the following: \\
$A_i = i \text{ for } i = 1,\ldots,t$ \\
$B_j = t + j \text{ for } j = 1,\ldots,t$ \\
$C_{ij} = 2t + (i-1)t + j \text{ for } i,j = 1,\ldots,t$ \\
$S = t^2 + 2t$ \\
$T = S + 2d + 1$

An instance of $2RmD$ with $t^2 + t$ resources and travel times between demands equal to zero  can be built as follows: resources $1,\ldots,t$ are of type 1 and the remaining $t^2$ resources are of type $0$. Let demands be of the form $(u_j,v_j)$, where $u_j$ is the start time, $v_j$ is the end time, and $v_j - u_j$ is the service time of demand $j$. Define the reward for demand $j$ to be the length of service times the number of resource types required.

The following demands require resource type $0$ and $1$:
$$(0,A_i) \text{ for } i = 1,\ldots,t$$
$$(A_i,C_{ij}) \text{ for } i,j = 1,\ldots,t$$
$$(S+d-c_k,T) \text{ for } k = 1,\ldots,t$$

The following demands require resource type $0$: 
$$t-1 \text{ times } (0,B_j) \text{ for } j = 1,\ldots,t$$
$$(B_j,C_{ij}) \text{ for } i,j = 1,\ldots,t$$
$$(C_{ij},S+a_i+b_j) \text{ for } i,j = 1,\ldots,t$$
$$(S+a_i+b_j,T-1) \text{ for } i,j = 1,\ldots,t$$
$$t^2-t \text{ times } (T-1,T)$$

The following demands require resource type 1:
$$(C_{ij},S+a_i+b_j) \text{ for } i,j = 1,\ldots,t$$

 We now show that there exists a feasible solution to N3DM if and only if all resources are serving demands during the interval [0,T] and there is no idle time. 
 
 Suppose that there exists a feasible schedule for a subset of demands, such that all resources are busy serving demands during the interval [0,T]. Demands $(0,A_i)$ must be scheduled to the first $t$ resources (since there are $t$ of these trips and they each require resource type 1) and resources $t+1,\ldots,2t$ of type $0$. These demands must be followed by some $(A_i,C_{ij})$ demands. Similarly, the $t^2-t$ demands $(0,B_j)$ must be scheduled to the remaining $t^2-t$ resources of type $0$, followed by some $(B_j,C_{ij})$ demands.

The first $t$ resources of type 1 and the next $t+1,\ldots,2t$ resources of type $0$ have the following schedules:
$$(0,A_i)(A_i,C_{ij})(C_{ij},S+a_i+b_j)(S+d-c_k,T)$$
where each $i$, $1\leq i \leq t$ occurs exactly once.

The remaining $t^2-t$ resources of type $0$ have schedules of the form
$$(0,B_j)(B_j,C_{ij})(C_{ij},S+a_i+b_j)(S+a_i+b_j,T-1)(T-1,T)$$
where each $j$, $1 \leq j \leq t$ occurs exactly $t-1$ times. Thus, among the demands $(C_{ij},S+a_i+b_j)$ that are served by the first $t$ resources, each $i$ and each $j$ occurs exactly once.

From the schedules of the first $t$ resources, we have $S+a_i+b_j = S+d-c_k$, i.e., $a_i + b_j + c_k = d$. So we define $\rho(i) = j$ and $\sigma(i) = k$ whenever demand $(C_{ij},S+a_i+b_j)$ is served followed by $(S+d-c_k,T)$ by resource $i, i = 1,\ldots,t$.

Conversely, given a feasible solution to N3DM, a feasible solution for this instance of $2RmD$ can be found, which is clearly optimal, since the resources are serving demand throughout the interval [0,T] with no idle time. \Halmos

\subsection{Proof of Lemma \ref{IPMR_MCF}}
\label{lemma_ipmr}
We first show that IPM-R can be modeled and solved as a MCF with node capacities. Define $G = (N,A,l,\mu,c,b,v)$ as follows, where $l,\mu$ are the lower and upper bound functions on arc capacity, $c$ is the arc cost function, $b$ is the node supplies function, and $v$ is the node capacities function:\\
\textbf{Nodes} \\
$N \coloneqq D^r \cup S^r \cup \{s^*,t\}$ \\
\textbf{Arcs:}\\
$(s^*,t)$ at cost 0 with capacity $(0,\infty)$ \\
$(s^*,s)$ at cost 0 with capacity $(0,l_s^r)$ $\forall s \in S^r$ \\
$(d,t)$ at cost 0 with capacity $(0,1)$ $\forall d \in D^r$ \\
$(s,d)$ if $B_{sd} = 1$ at cost 0 with capacity $(0,1)$ $\forall s \in S^r, d \in D^r$ \\
$(i,j)$ if $A_{ij} = 1$ at cost 0 with capacity $(0,1)$ $\forall i,j \in D^r$ \\
\textbf{Supplies (b):} \\
$b(i) = 0$ $\forall i \in S^r \cup D^r$ \\
$-b(t) = b(s^*) = \sum_{s \in S^r}l_s^r$ \\
\textbf{Node Capacities (v):} \\
$v(i) = y_i \text{ if } i \in D^r$

To convert this MCF with node capacities into a general MCF, we can perform the following transformation (\cite{ciupala}). We redefine $G$ to be a MCF $G' = (N',A',l',\mu',c',b')$. \\
\textbf{Nodes:} \\
$N' = N_1 \cup N_2$ where $N_1 \coloneqq \{i'|i \in N\}$ and $N_2 \coloneqq \{i''|i \in D^r\}$. \\
\textbf{Arcs:} \\
$A' = A_1 \cup A_2 \cup A_3 \cup A_4 \cup A_5$ where \\
$A_1 \coloneqq \{(i',j')|(i,j) \in A, i,j \notin D^r\}$ \\
$A_2 \coloneqq \{(j',i')|(j,i) \in A, j \notin D^r, i \in D^r\}$ \\
$A_3 \coloneqq \{(i'',j')|(i,j) \in A, i,j \in D^r\}$ \\
$A_4 \coloneqq \{(i'',j')|(i,j) \in A, i \in D^r, j \notin D^r\}$ \\
$A_5 \coloneqq \{(i',i'')| i \in D^r\}$ \\
Note that for any arc $(i,j) \in A_1 \cup A_2 \cup A_3 \cup A_4$, the cost and capacities remain the same as previously defined for the corresponding arcs in $A$. The new arcs $(i',i'') \in A_5$ have cost 0 and capacity $(v(i),v(i))$.\\
\textbf{Supplies (b'):} \\
$b'(i') = b(i)$ $\forall i' \in N_1$ \\
$b'(i'') = 0$ $\forall i'' \in N_2$ \\
Then we have the following MCF formulation:
\begin{align*}
         \min & \hspace{1em} 0&
\end{align*}
s.t.
\begin{align*}
    \sum_{j:(i,j) \in A'}x_{ij}^r - \sum_{j:(j,i) \in A'}x_{ji}^r = b'(i) &  & i \in N' \\
    l'(i,j) \leq x_{ij}^r \leq \mu'(i,j) &  & (i,j) \in A'
\end{align*}
Thus, we have shown that IPM-R can be modeled and solved as a MCF. \Halmos

\subsection{Proof of Theorem \ref{non_approx}}
\label{pf_non_approx}
We prove that $mRmD$ cannot be approximated within $O(\min\{|R|^{\frac{1}{2}-\epsilon},|D|^{1-\epsilon}\})$, where $|R|$ is the number of resource types and $|D|$ is the number of demands, by reduction using Maximum Set Packing.

\noindent \textbf{Instance of Maximum Set Packing} \\
Given a universe $U$, a family $S$ of subsets of $U$, and an integer $n$, is there a subfamily $C \subseteq S$ of sets such that all sets in $C$ are pairwise disjoint and $|C| \geq n$?

Maximum Set Packing is well known to be NP-Complete and cannot be approximated within $O(N^{1-\epsilon})$ unless $NP \subseteq ZPP$ (\cite{hastad}). Given an instance of Maximum Set Packing, we can create an instance of $mRmD$ with $|U|$ resources and $|S|$ demands in polynomial time as follows:

\noindent Define $r_u$ as a resource of type $u$ $\forall u \in U$. \\
Define $d_s$ as a demand that requires resource type $i$ if $i \in s$ $\forall s \in S$. \\
Let all demands and resources be at the same location (i.e., zero travel times) and all demand intervals be equivalent. Further, set all demand rewards equal to 1.

We now show that there exists a feasible solution to Maximum Set Packing if and only if there exists a feasible solution to our problem with an objective value greater than or equal to $n$.

Suppose that we have a feasible solution to Maximum Set Packing. That is, $\exists C \subseteq S$ such that $|C| \geq n$ and all sets in $C$ are pairwise disjoint. Then, for each $x \in C$, $d_x$ can be served by our resources since they are disjoint in the resources that they require. Therefore, we can serve at least $n$ demands and since all rewards are equal to 1, our objective must be greater than or equal to $n$.

Conversely, suppose we have a feasible solution to our problem such that the objective is greater than or equal to $n$. Since all rewards are equal to 1, this implies that we have met at least $n$ demands. Further, since we only have one resource of each type, these $n$ demands must have disjoint resource requirements. For each demand $d_x$ that was served, the corresponding set of resource requirements $x \in S$ can be added to $C$. Therefore, we have created a subfamily $C \subseteq S$ of disjoint sets where $|C| \geq n$ and so Maximum Set Packing is feasible.

Since the transformation of Maximum Set Packing to an instance of $mRmD$ preserved the objective value, then approximation results are also preserved. That is, $mRmD$ cannot be approximated within $O(\min\{|R|^{\frac{1}{2}-\epsilon},|D|^{1-\epsilon}\})$ unless $NP \subseteq ZPP$ (for every $\epsilon > 0$), where $|R|$ is the number of resource types and $|D|$ is the number of demands.  \Halmos

\subsection{Proof of Theorem \ref{mrt-1or2}}
\label{pf_mrt-1or2}
We prove that the special case of $mR\{\text{1 or 2}\}D$ where there are an arbitrary number of resource types $R$, a special resource $r'\in R$, and each demand requires either $r',r$ or $ \{r',r\}, r \in R - r'$ is NP-hard by reduction using N3DM. The construction is similar to the proof in \cite{aircraft}. Given an instance of N3DM, we define the following: \\
$A_i = i \text{ for } i = 1,\ldots,t$ \\
$B_j = t + j \text{ for } j = 1,\ldots,t$ \\
$C_{ij} = 2t + (i-1)t + j \text{ for } i,j = 1,\ldots,t$ \\
$S = t^2 + 2t$ \\
$T = S + 2d + t + 1$

An instance of the special case of $mR\{\text{1 or 2}\}D$ with $t^2 + t$ resources and travel times between demands equal to zero can be built as follows: resources $1,\ldots,t$ are of type $i, i = 1, \ldots, t$ and the remaining $t^2$ resources are of type $t'$. Let demands be of the form $(u_j,v_j)$, where $u_j$ is the start time, $v_j$ is the end time, and $v_j - u_j$ is the service time of demand $j$. Define the reward for demand $j$ to be the length of service times the number of resource types required.

The following demands require resource type $t'$ and $i$:
$$(0,A_i) \text{ for } i = 1,\ldots,t$$
$$(A_i,C_{ij}) \text{ for } i,j = 1,\ldots,t$$
$$(S+d-c_k,T-k-1) \text{ for } i,k = 1,\ldots,t$$

The following demands require resource type $t'$: 
$$t-1 \text{ times } (0,B_j) \text{ for } j = 1,\ldots,t$$
$$(B_j,C_{ij}) \text{ for } i,j = 1,\ldots,t$$
$$(C_{ij},S+a_i+b_j) \text{ for } i,j = 1,\ldots,t$$
$$(S+a_i+b_j,T-1) \text{ for } i,j = 1,\ldots,t$$
$$t^2-t \text{ times } (T-1,T)$$
$$(T-k-1,T) \text{ for } k = 1,\ldots,t$$

The following demands require resource type $i$:
$$(C_{ij},S+a_i+b_j) \text{ for } i,j = 1,\ldots,t$$
$$(T-k-1,T) \text{ for } i,k = 1,\ldots,t$$

 We now show that there exists a feasible solution to N3DM if and only if all resources are serving demands during the interval [0,T] and there is no idle time. 

Suppose that there exists a feasible schedule for a subset of demands, such that all resources are busy serving demands during the interval [0,T]. Demands $(0,A_i)$ must be scheduled to the first $t$ resources (since there are $t$ of these trips and they each require resource type $i, i = 1,\ldots t$) and resources $t+1,\ldots,2t$ of type $t'$. These demands must be followed by some $(A_i,C_{ij})$ demands. Similarly, the $t^2-t$ demands $(0,B_j)$ must be scheduled to the remaining $t^2-t$ resources of type $t'$, followed by some $(B_j,C_{ij})$ demands.

The first $t$ resources of type $i, i = 1,\ldots,t$ and the next $t+1,\ldots,2t$ resources of type $t'$ have the following schedules:
$$(0,A_i)(A_i,C_{ij})(C_{ij},S+a_i+b_j)(S+d-c_k,T-k-1)(T-k-1,T)$$
where each $i$ and each $k$, $1\leq i,k \leq t$ occur exactly once.

The remaining $t^2-t$ resources of type $t'$ have schedules of the form
$$(0,B_j)(B_j,C_{ij})(C_{ij},S+a_i+b_j)(S+a_i+b_j,T-1)(T-1,T)$$
where each $j$, $1 \leq j \leq t$ occurs exactly $t-1$ times. Thus, among the demands $(C_{ij},S+a_i+b_j)$ that are served by the first $t$ resources, each $i,j,$ and $k$ occurs exactly once.

From the schedules of the first $t$ resources, we have $S+a_i+b_j = S+d-c_k$, i.e., $a_i + b_j + c_k = d$. So we define $\rho(i) = j$ and $\sigma(i) = k$ whenever demand $(C_{ij},S+a_i+b_j)$ is served followed by $(S+d-c_k,T-k-1)$ by resource $i, i = 1,\ldots,t$.

Conversely, given a feasible solution to N3DM, a feasible solution for this instance of the special case of $mR\{\text{1 or 2}\}D$ can be found, which is clearly optimal, since the resources are serving demand throughout the interval [0,T] with no idle time. \Halmos

\subsection{Proof of Theorem \ref{m_1orall}}
\label{pf_m_1orall}
For an instance $\mathcal{I}$ of the special case of $mR\{\text{1 or all}\}D$, let $m_i$ be the number of resources of type $i$ for $i = 1, \ldots, |R|$. Since there is zero travel time between demands (i.e., all demands are at the same location), all demands have the same service start time, and all service times are infinite, then each resource can only be assigned a single job. 

First, divide the list of demands $D$ into subsets $D_0,D_1,\ldots,D_{|R|}$ where subset $D_0$ is the demands that require all resource types and subsets $D_i$ for $i = 1, \ldots, |R|$ are the demands that only require resource type $i$. Order the demands in each subset by decreasing reward.  We can find the optimal objective value $OPT(\mathcal{I})$ and its corresponding solution by Algorithm \alge $ $, which runs polynomially in number of resources. Note the solution is constructed by serving all demands whose reward was added to the objective value $OPT(\mathcal{I})$. \Halmos

\begin{algorithm}[!htb]
\renewcommand{\thealgocf}{}
\caption{\alge}
\label{algo_notravel_inf}
 \SetKwInOut{Input}{Input}
 \SetKwInput{Output}{Output}
 \Input{$\mathcal{I}$, $m_i$, $i = 1,2,\ldots,|R|$}
$m \leftarrow \min\{m_i|i\in\{1,2,\ldots,|R|\}\}$ \\
$Opt(\mathcal{I}) \leftarrow 0$ \\
  \For{$1 \leq i \leq |R|$}{
    \If{$m < m_i$}{
        Increase $OPT(\mathcal{I})$ by rewards of first $m_i-m$ entries of $D_i$. \\
        Remove first $m_i - m$ entries from $D_i$.
    }
  }
  \For{$1 \leq i \leq m$}{
    Let $r(D_i[0])=$ reward of first entry in $D_i$ for $i = 0,1,\ldots,|R|$.\\
    \If{$\sum_{i=1}^{|R|} r(D_i[0]) \leq r(D_0[0])$}{
        Increase $OPT(\mathcal{I})$ by $r(D_0[0])$. \\
        Remove $D_0[0]$ from $D_0$. 
    }
    \Else{
        Increase $OPT(\mathcal{I})$ by $r(D_i[0])$ for $i = 1,2,\ldots,|R|$. \\
        Remove $D_i[0]$ from $D_i$ for $i=1,2,\ldots, |R|$.
    }
  }
  \Output{$OPT(\mathcal{I})$}
\end{algorithm}

\subsection{Proof of Theorem \ref{1rt-1d}}
\label{pf_1rt-1d}
We prove that $1R1D$ is solvable in polynomial time. Consider $m$ resources and $n$ demands. Without loss of generality, assume that each resource has its own starting location and are labeled by resource number (i.e., $s = 1, \ldots, m$). This problem can be modeled as a MCF on a directed acyclic graph with the following nodes and arcs: \\
\textbf{Nodes:}\\
For each resource there is a node $b_i$, $i = 1, \ldots m$. \\
For each demand there are two nodes $u_j$ and $v_j$, $j = 1, \ldots n$. \\
There is a dummy source node $s$ and a dummy sink node $t$. \\
\textbf{Arcs:} \\
$(s,b_i)$ and $(b_i,t)$ at cost $0$, $i=1, \ldots, m$. \\
$(v_j,t)$ at cost $0$, $j = 1, \ldots n$. \\
$(b_i,u_j)$ if $B_{ij}=1$ at cost 0, $i=1,\ldots,m, j=1,\ldots,n$. \\
$(v_j,u_k)$ if $A_{jk}=1$ at cost 0, $j,k=1,\ldots,n$. \\
$(u_j,v_j)$ at cost $-w_j$, $j=1,\ldots,n$. \\
All arcs have capacity 1. Let the supply at node $s$ be $m$, supply at node $t$ be $-m$ and supply for all other nodes be zero. Further, note that since the graph is directed and acyclic, we can transform the arc costs to be non-negative. Because of the capacities of the arcs, each arc will have either no flow or flow of 1 unit in the optimal solution. It is clear that there is a one-to-one correspondence between the MCF from $s$ to $t$ on this constructed graph and the optimum allocation for a single type of resource. If there is a positive flow on the arc $(b_i,t)$, this means that resource $i$ was not used to meet demand. If there is a positive flow on the arc $(b_i,u_j)$, this means that demand $j$ is the first demand met by resource $i$. Positive flow on the arc $(v_j,u_k)$ means that demand $k$ is met immediately after demand $j$ by the same resource. Finally, positive flow on the arc $(u_j,v_j)$ means that demand $j$ was met. \Halmos

\subsection{Observation 1}
\label{obs_1}
Note that if travel times are zero, then $1R1D$ can be formulated as a maximum weighted coloring problem on an interval graph. Consider $m$ resources and $n$ demands. We can now create an (undirected) interval graph $G = (N,A)$ as follows: \\
\textbf{Nodes:} \\
For each demand, create a node $d_i, i = 1, \ldots, n$. \\
\textbf{Arcs:} \\
$(d_i,d_j)$ if $\tau_i + \Delta_i > \tau_j \hspace{1em} \forall i,j \in \{1,2, \ldots, n\}$ \\
Then, we can solve the maximum weighted $m$-coloring problem on interval graph $G$ using the following integer program:
\begin{align*}
    \max  \sum_{i=1}^{n} w_id_i &
\end{align*}
s.t.
\begin{align*}
\sum_{i \in S}d_i \leq m & & \forall  S \text{ s.t. $S$ is a maximal clique of size $\geq m$} \\
    d_i \in \{0,1\} & &\forall  i = 1,\ldots, n
\end{align*}
Since the nodes versus cliques matrix of an interval graph is totally unimodular (\cite{mertzios}), then we can relax the above formulation to a linear program and still maintain integral optimal solutions. 

\section{Computational Study}
\subsection{Bicriteria Results}
\label{bicriteria_ext}
Figures \ref{bicriteria_3R} and \ref{bicriteria_2R} present bicriteria results to compare the trade-off between increasing the number of resources and its impact on the objective function value (i.e., demands met). These results could provide insights for decision-makers when determining resource capacity and desired levels of demand satisfaction. In Figure \ref{bicriteria_3R}, we examine the $mRmD$ problem with 3 resource types, 400 demands, and 30 resources. Each line in the graph represents an instance of this problem. Starting at 10 resources of each type, we increase this amount by 0 to 15 and record the objective function value. As shown, there is a diminishing return as resources are increased. On average, after increasing the resources by about 6 or 7, we reach the maximum objective function value for each instance studied. Figure \ref{bicriteria_2R} looks at the $2RmD$ problem with 200 demands and 14 resources. In this case, we increase the number of resources of each type separately by 0 to 11 and present the resulting objective function value. Again, we see a diminishing return as resources are increased. It is worth noting that jointly increasing resources (i.e., increasing both resource types) is more profitable than only increasing one resource type. For example, increasing resource type 1 and type 2 by 3 and 4 resources, respectively, produced an objective of 13606 whereas only increasing resource type 1 by 8 led to an objective of 12844.

\begin{figure}
\centering
\begin{minipage}{.5\textwidth}
  \centering
   \caption{Bicriteria results for $mRmD$}
  \includegraphics[width=\linewidth]{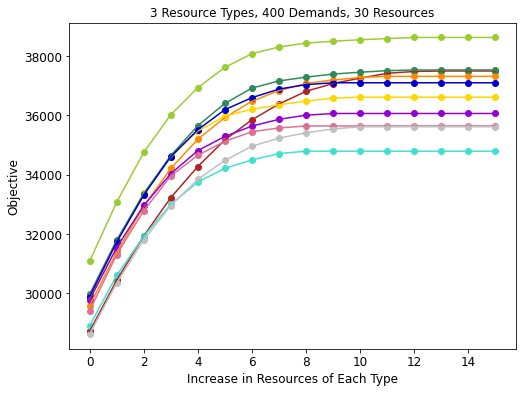}
  \label{bicriteria_3R}
\end{minipage}%
\begin{minipage}{.5\textwidth}
  \centering
   \caption{Bicriteria results for $2RmD$}
  \includegraphics[width=\linewidth]{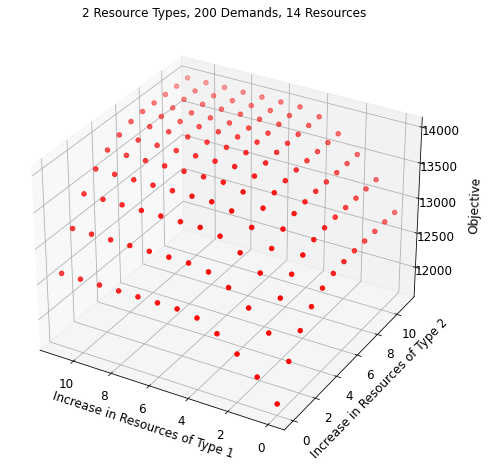}
  \label{bicriteria_2R}
\end{minipage}
\end{figure}

\subsection{Travel Costs}
\label{comp_ext}
Let $BFlow$ and $CFlow$ represent solutions in which the flow variables ($x$) are restricted to be binary and those in which the flow variables are relaxed to be continuous, respectively.
Tables \ref{ext_results_1} and \ref{ext_results_2} present a comparison of the run times (in seconds) for $BFlow$ and $CFlow$ under all instances considered, where $|R|$ is the number of resource types, $|D|$ is the number of demands, and $L$ is the total number of resources. Note that \textit{Scaled Demands} refers to instances where reward for demand was multiplied by a factor of 100.

Computational experiments show that for small and medium sized problems (less than 5 resource types), IPM-C can be efficiently solved by Gurobi and the optimal solution can be found in reasonable time. As the size of the problem increases (i.e., greater demands or larger number of resource types), the instances become harder to solve. For 2-4 resource types, on average, $BFlow$ performs better than $CFlow$ when demands are not scaled, whereas when demands are scaled, CFlow performs better than BFlow. For 5-7 resource types, on average, $CFlow$ performs better than $BFlow$ in both scaled and unscaled demand scenarios. However, for some instances, $BFlow$ is drastically better than $CFlow$ (i.e., 7 resource types and 700/800 demands for unscaled/scaled demands, respectively). In general instances with scaled demands perform better than those in which demands are not scaled; this could be a result of less importance given to travel costs in the objective.

\begin{table}[!htb]
    \begin{minipage}{0.5\textwidth}
        \centering
        \caption{Results for number of resource types equal to 2,3 and 4}
        \label{ext_results_1}
        \resizebox{!}{.15\paperheight}{%
\begin{tabular}{|c|c|c|c|c|c|c|}
\hline
\textbf{$|R|$} & \textbf{$|D|$} & \textbf{$L$} & \textbf{\begin{tabular}[c]{@{}c@{}}Run time\\ (BFlow)\end{tabular}} & \textbf{\begin{tabular}[c]{@{}c@{}}Run time\\ (CFlow)\end{tabular}} & \textbf{\begin{tabular}[c]{@{}c@{}}Run time\\ (BFlow)\\ (Scaled Demands)\end{tabular}} & \textbf{\begin{tabular}[c]{@{}c@{}}Run time\\ (CFlow)\\ (Scaled Demands)\end{tabular}} \\ \hline
2 & 100 & 10 & 0.11 & 0.09 & 0.16 & 0.10 \\ \hline
2 & 200 & 14 & 0.28 & 0.25 & 0.33 & 0.31 \\ \hline
2 & 300 & 18 & 0.71 & 0.72 & 1.10 & 0.69 \\ \hline
2 & 400 & 22 & 1.54 & 2.01 & 2.24 & 1.70 \\ \hline
2 & 500 & 26 & 3.47 & 3.96 & 3.84 & 3.47 \\ \hline
2 & 600 & 30 & 5.49 & 6.90 & 6.79 & 6.89 \\ \hline
2 & 700 & 34 & 8.35 & 13.77 & 8.05 & 9.80 \\ \hline
2 & 800 & 38 & 12.69 & 19.09 & 15.16 & 14.44 \\ \hline
3 & 100 & 12 & 0.09 & 0.09 & 0.15 & 0.10 \\ \hline
3 & 200 & 18 & 0.49 & 0.39 & 0.57 & 0.59 \\ \hline
3 & 300 & 24 & 1.01 & 1.15 & 1.48 & 1.10 \\ \hline
3 & 400 & 30 & 3.81 & 3.63 & 3.78 & 3.42 \\ \hline
3 & 500 & 36 & 7.29 & 8.23 & 6.72 & 5.90 \\ \hline
3 & 600 & 42 & 9.98 & 13.51 & 13.54 & 13.69 \\ \hline
3 & 700 & 48 & 19.40 & 30.17 & 20.65 & 21.99 \\ \hline
3 & 800 & 54 & 43.85 & 65.33 & 25.16 & 61.34 \\ \hline
4 & 100 & 16 & 0.28 & 0.14 & 0.24 & 0.13 \\ \hline
4 & 200 & 24 & 0.97 & 0.82 & 1.08 & 0.55 \\ \hline
4 & 300 & 32 & 2.63 & 2.19 & 2.29 & 1.56 \\ \hline
4 & 400 & 40 & 6.22 & 5.41 & 5.31 & 4.67 \\ \hline
4 & 500 & 48 & 19.49 & 15.28 & 15.39 & 14.26 \\ \hline
4 & 600 & 56 & 26.87 & 30.03 & 23.84 & 22.85 \\ \hline
4 & 700 & 64 & 64.84 & 96.78 & 40.41 & 38.23 \\ \hline
4 & 800 & 72 & 89.05 & 125.97 & 52.18 & 77.69 \\ \hline
\end{tabular}}
    \end{minipage}
    \hfill
    \begin{minipage}{0.5\textwidth}
        \centering
        \caption{Results for number of resource types equal to 5,6 and 7}
        \label{ext_results_2}
        \resizebox{!}{.15\paperheight}{%
       \begin{tabular}{|c|c|c|c|c|c|c|}
\hline
\textbf{$|R|$} & \textbf{$|D|$} & \textbf{$L$} & \textbf{\begin{tabular}[c]{@{}c@{}}Run time\\ (BFlow)\end{tabular}} & \textbf{\begin{tabular}[c]{@{}c@{}}Run time\\ (CFlow)\end{tabular}} & \textbf{\begin{tabular}[c]{@{}c@{}}Run time\\ (BFlow)\\ (Scaled Demands)\end{tabular}} & \textbf{\begin{tabular}[c]{@{}c@{}}Run time\\ (CFlow)\\ (Scaled Demands)\end{tabular}} \\ \hline
5 & 100 & 20 & 0.19 & 0.19 & 0.13 & 0.18 \\ \hline
5 & 200 & 30 & 2.00 & 1.01 & 1.10 & 0.80 \\ \hline
5 & 300 & 40 & 6.36 & 3.68 & 3.39 & 2.74 \\ \hline
5 & 400 & 50 & 14.79 & 9.20 & 9.04 & 7.76 \\ \hline
5 & 500 & 60 & 35.01 & 32.58 & 21.28 & 17.83 \\ \hline
5 & 600 & 70 & 50.66 & 70.35 & 51.90 & 73.11 \\ \hline
5 & 700 & 80 & 96.25 & 172.39 & 47.01 & 53.80 \\ \hline
5 & 800 & 90 & 453.05 & 1014.62 & 134.52 & 202.69 \\ \hline
6 & 100 & 24 & 0.44 & 0.22 & 0.29 & 0.26 \\ \hline
6 & 200 & 36 & 2.49 & 1.25 & 1.23 & 1.40 \\ \hline
6 & 300 & 48 & 7.08 & 5.43 & 5.44 & 4.44 \\ \hline
6 & 400 & 60 & 18.03 & 15.78 & 11.36 & 10.40 \\ \hline
6 & 500 & 72 & 53.19 & 55.05 & 30.33 & 33.52 \\ \hline
6 & 600 & 84 & 134.02 & 107.60 & 137.94 & 145.15 \\ \hline
6 & 700 & 96 & 195.10 & 237.03 & 114.26 & 181.08 \\ \hline
6 & 800 & 108 & 537.75 & 1039.68 & 232.16 & 285.46 \\ \hline
7 & 100 & 28 & 0.33 & 0.25 & 0.29 & 0.20 \\ \hline
7 & 200 & 42 & 2.06 & 1.63 & 1.85 & 1.27 \\ \hline
7 & 300 & 56 & 13.76 & 8.58 & 9.92 & 5.92 \\ \hline
7 & 400 & 70 & 125.48 & 55.14 & 23.67 & 17.23 \\ \hline
7 & 500 & 84 & 161.38 & 112.70 & 67.47 & 51.95 \\ \hline
7 & 600 & 98 & 421.19 & 389.39 & 133.33 & 167.05 \\ \hline
7 & 700 & 112 & 989.80 & 1629.80 & 306.81 & 279.56 \\ \hline
7 & 800 & 126 & 1485.61 & 1398.86 & 359.35 & 534.93 \\ \hline
\end{tabular}}
    \end{minipage}
\end{table}

\section{Extension 2}
\label{extension}
We present an extension to $mRmD$. In this extension, resource types have both an origin and destination location, i.e., resources may not go directly from demand incident to demand incident. For an overview of the additional notation needed for this formulation, please refer to Table \ref{formulation_ext}.

In order to construct feasible schedules for resources, we create 0-1 matrices for each resource type $r \in R$, $A^r$, and a 0-1 matrix, $B$. Matrix $A^r$ has a row and column for each demand incident and $A^r_{ij} = 1$ if and only if demand incident $j$ can be served after demand incident $i$ by resource type $r$ and demands $i$ and $j$ both require resource type $r$. Matrix $B$ has a row for each resource starting location and a column for each demand incident and $B_{sd}=1$ if and only if a resource starting at location $s$ can serve demand $d$ first. That is, we have the following pre-processing matrices:\begin{align*}
    A^r_{ij} &= 
    \begin{cases}
    1  \text{ iff} \hspace{1em} \tau_i + \Delta_i + f(a_i^r,b_i^r) + f(b_i^r,a_j^r) \leq \tau_j, \hspace{1em} r \in M_i \cap M_j  \\
    0 \text{ otherwise}
    \end{cases}
    &  & i,j \in D, r \in R \\
    B_{sd} &=
    \begin{cases}
    1 \text{ iff} \hspace{1em}  f(s,a_d^r) \leq \tau_d, s \in S^r,r \in M_d \\
    0 \text{ otherwise}
    \end{cases}
    &  & s \in S^r, r \in R, d \in D
\end{align*}  
The formulation is the same as IPM, where variable construction relies on matrices $A^r$ and $B$.
\begin{table}[!htb]
    \centering
    \small\setlength\tabcolsep{4.5pt}
    \caption{Additional notation for Extension 2}
    \begin{tabular}{l l}
    \hline 
    \hline
    $V$ & Set of nodes \\
    $R$ & Set of resource types \\
    $f(u,v)$ & Travel time between $u,v \in V$ \\
    $(a_d^r,b_d^r)$ & Origin-destination pair for resource type $r \in M_d$ needed by $d \in D$ \\
    \hline
    \hline
    \end{tabular}
    \label{formulation_ext}
\end{table}

\end{APPENDICES}


\bibliographystyle{informs2014} 
\bibliography{bibliography.bib} 


\end{document}